\documentclass[12pt]{article}
\usepackage{amssymb}
\usepackage{amsmath}
\usepackage{multirow}
\usepackage{hyperref}
\usepackage{hyperref}
\hypersetup{
  colorlinks   = true, 
  urlcolor     = blue, 
  linkcolor    = blue, 
  citecolor   = red 
}

\newtheorem{theorem}{Theorem}[section]

\newtheorem{lemma}[theorem]{Lemma}

\begin{document}

\title{A Jacobi theta series and its transformation laws}
\author{Matthew Krauel
}

\date{}
\maketitle

\abstract 
\noindent
We consider a generalization of Jacobi theta series and show that every such function is a quasi-Jacobi form. Under certain conditions we establish transformation laws for these functions with respect to the Jacobi group and prove such functions are Jacobi forms. In establishing these results we construct other functions which are also Jacobi forms. These results are motivated by applications in the theory of vertex operator algebras.

\section{Introduction}

 \indent  Let $Q$ be a positive definite integral quadratic form, and $B$ be the associated bilinear form so that $2Q(x)=B(x,x)$. Fix $h\in \mathbb{Z}^f$ and let $A$ be the matrix of $Q$ with even rank $f=2r$. It is well known (page $81$ of \cite{EZ}, for example) that the \textit{Jacobi theta series} $\sum_{m\in \mathbb{Z}^f} q^{Q(m)} \zeta^{B(m,h)}$ is a Jacobi form of weight $r$ and index $Q(h)$ on the full Jacobi group $J^1 =\text{SL}_2 (\mathbb{Z})\ltimes (\mathbb{Z}\times \mathbb{Z})$. Here, $\zeta =e^{2\pi i z}$ with $z\in \mathbb{C}$, and $q=e^{2\pi i \tau}$ with $\tau$ in the complex upper half-plane $\mathbb{H}$. \\
 
 \indent Fix an element $v\in \mathbb{C}^f$. Let $\underline{z}=(z_1 ,\dots ,z_n) \in \mathbb{C}^n$, and set $\underline{h}=(h_1 ,\dots ,h_n)$ for linearly independent elements $h_1 ,\dots ,h_n \in \mathbb{Z}^f$. In this paper we consider functions of the form
 \begin{equation}
 \theta_{\underline{h}} (Q,v,k,\tau ,\underline{z}):=\sum_{m\in \mathbb{Z}^f} B(v,m)^k q^{Q(m)} \zeta_1^{B(m,h_1)} \cdots \zeta_n^{B(m,h_n)}, \label{thetaspehricals}
 \end{equation}
 for any natural number $k$. We show such functions are quasi-Jacobi forms of weight $k+r$ for all $v \in \mathbb{C}^f$. In the case $v$ is restricted to $\mathbb{R}^f$ and satisfies $Q(v)=0$, the positive definiteness of $Q$ implies $v=0$. However, there may be nonzero $v\in \mathbb{C}^f$ such that $Q(v)=0$. For $v\in \mathbb{C}^f$ such that $Q(v)=0$ and $B(v,h_j)=0$ for all $j$, we develop transformation laws with respect to a subgroup of the Jacobi group $J^n =\text{SL}_2 (\mathbb{Z})\ltimes (\mathbb{Z}^n \times \mathbb{Z}^n)$ and show these functions are Jacobi forms on this subgroup. Precise definitions of Jacobi forms and quasi-Jacobi forms are given in Section \ref{Preliminaries}.\\
 \indent In the special case $\underline{z}=0$, such functions are considered in \cite{DM-Theta} and transformation laws with respect to $\text{SL}_2 (\mathbb{Z})$ are developed which expands on work of \cite{Hecke-Positiven,Ogg-Modular,Schoeneberg-Das,Schoeneberg-Elliptic} for similar functions. Many of the techniques used in the proofs below are attributed to these authors. We will reference them often, and attempt to maintain the notation developed in these works (especially \cite{DM-Theta}). Other techniques have been developed to establish transformation laws for theta series of higher degree \cite{AM-BehaviorI}.\\
 \indent Holomorphic Jacobi forms of higher degree have been considered by Ziegler in \cite{Ziegler-JacobiForms}, where Jacobi theta series consisting of two complex matrix variables are shown to be examples. Richter generalizes some of this work in \cite{Richter-OnTransformations}. However, the literature lacks such results for Jacobi theta series as defined in (\ref{thetaspehricals}). This paper fills this gap and is motivated by the occurrence of functions of the form (\ref{thetaspehricals}) and (\ref{PsiFunctions}) in the theory of vertex operator algebras in work similar to those of \cite{DMN-quasi} and \cite{KrauelMasonI}. This work will appear elsewhere. \\

 Let $\Gamma :=\text{SL}_2 (\mathbb{Z})$ and $\Gamma_0 (N)$ be the subgroup of $\Gamma$ defined by
 \begin{equation*}
 \Gamma_0 (N)=\left \{ \begin{pmatrix} a&b \\c&d \end{pmatrix} \in \Gamma \mid c \equiv 0 \mod N \right \}. \label{Gamma0(N)}
 \end{equation*}
 We denote the subgroup $\Gamma_0 (N)\ltimes (\mathbb{Z}^n \times \mathbb{Z}^n)$ of $J^n$ by $J^n_0 (N)$ and the spaces of Jacobi forms and quasi-Jacobi forms of weight $k$ and index $F$ on the group $J^n_0 (N)$ by $\mathcal{J}^{n,k,F}_0 (N)$ and $\mathcal{Q}^{n,k,F}_0 (N)$, respectively. Let $G$ be the $n\times n$ matrix whose $(i,j)$-th entry is given by $B(h_i,h_j)$. In the case $n=f$, $h_1 ,\dots ,h_n$ is a basis of $\mathbb{Z}^f$ and $G$ is the Gram matrix associated to this set of elements.\\
 	\indent For a column vector $\underline{\alpha}=(\alpha_1 ,\dots ,\alpha_n)$, set $G[\underline{\alpha}]:=\underline{\alpha}^T G\underline{\alpha}$ where $x^T$ denotes the transpose of $x$. We define $\epsilon (d)$ for $d>0$ by $\epsilon (d)=\left( \frac{(-1)^r \det (A)}{d}\right)$, and $\epsilon (-d)=(-1)^r \epsilon (d)$. The function $\epsilon$ is a Dirichlet character (see page $216$ in \cite{Schoeneberg-Elliptic} for example). For a matrix $\gamma =\left( \begin{smallmatrix} a& b\\c &d \end{smallmatrix} \right)$, we often write $\gamma \tau$ and $\gamma \underline{z}$ to denote $\frac{a\tau +b}{c\tau +d}$ and $\frac{\underline{z}}{c\tau +d}$, respectively.\\
 \indent We are now in position to state our first result.
\begin{theorem}\label{theorem1}
For any $v$, $\underline{h}$, and $k$ as defined above, we have
\begin{equation}
\theta_{\underline{h}}(Q,v,k,\tau ,\underline{z}) \in \mathcal{Q}^{n,k+r,F}_0 (N) \label{thm1a}
\end{equation}
and $\theta_{\underline{h}}$ has a Fourier expansion of the form
\begin{equation}
 \theta_{\underline{h}}(Q,v,k,\tau ,\underline{z})=\underset{4\ell -(G/2)^{-1}[\underline{\nu}] \geq 0}{\sum_{\underline{\nu}\in \mathbb{Z}^n, \ell \in \mathbb{Q},}} c(\ell ,\underline{\nu}) q^\ell \exp \left( 2\pi i (\underline{z}^T \underline{\nu}) \right), \label{Expansion}
\end{equation}
where $\ell \geq 0$ and $c(\ell ,\underline{\nu})$ are scalars. In particular, $\ell =Q(m)\in \mathbb{Z}$ and $\underline{\nu}=(B(m,h_1), \dots ,B(m,h_n))$.\\
 \indent If in addition $Q(v)=0$ and $B(v,h_j)=0$ for all $1\leq j\leq n$, then for any $\gamma =\left( \begin{smallmatrix} a& b\\c &d \end{smallmatrix} \right) \in \Gamma_0 (N)$ and $(\underline{\lambda} ,\underline{\mu}) \in \mathbb{Z}^n \times \mathbb{Z}^n$ we have
 \begin{equation}
 \theta_{\underline{h}} \left(Q,v,k,\gamma \tau , \gamma \underline{z} \right)
 =\epsilon (d) \exp \left(\frac{\pi icG[\underline{z}]}{c\tau +d} \right) (c\tau +d)^{k+r} \theta_{\underline{h}} (Q,v,k,\tau ,\underline{z}), \label{ThetaTrans1}
 \end{equation}
  and
 \begin{equation}
 \theta_{\underline{h}} (Q,v,k,\tau ,\underline{z}+\underline{\lambda} \tau +\underline{\mu})=\exp \left(-\pi i(G[\underline{\lambda}] \tau +2\underline{z}^tT\underline{\lambda}) \right)\theta_{\underline{h}} (Q,v,k,\tau ,\underline{z}). \label{ThetaTrans2}
 \end{equation}
\end{theorem}

 \indent The  expansion (\ref{Expansion}) along with the transformation laws (\ref{ThetaTrans1}) and (\ref{ThetaTrans2}) imply that if $Q(v)=0$ and $B(v,h_j)=0$ for all $1\leq j\leq n$, then (\ref{thetaspehricals}) is a Jacobi form of weight $k+r$, index $G/2$, and character $\epsilon$. Theorem \ref{theorem1} is proved in Section \ref{ProofTheoremsMajor}. \\
 
 \indent In proving Theorem \ref{theorem1} we consider functions of the form
 \begin{equation}
 \Psi_{\underline{h}}(Q,v,k,\tau ,\underline{z}) :=\sum_{t=0}^{\lfloor k/2 \rfloor} \delta (t,k) (2Q(v)E_2 (\tau))^t \theta_{\underline{h}}(Q,v,k-2t,\tau ,\underline{z}), \label{PsiFunctions}
 \end{equation}
 where $E_2 (\tau)$ is the usual modular Eisenstein series of weight $2$, $\lfloor k/2 \rfloor$ denotes the greatest integer less than or equal to $k/2$, and $\delta (t,k)$ is defined by
 \begin{equation}
 \delta (t,k)= \frac{k!}{2^t t!(k-2t)!}. \label{rhoscalar}
 \end{equation} 
 We establish the following theorem in Section \ref{ProofTheoremsMajor}.
 
\begin{theorem} \label{theorem2}
Suppose $B(v,h_j)=0$ for all $1\leq j\leq n$. Then for any $\gamma =\left( \begin{smallmatrix} a& b\\c &d \end{smallmatrix} \right) \in \Gamma_0 (N)$ and $[\underline{\lambda} ,\underline{\mu}]\in (\mathbb{Z}^n \times \mathbb{Z}^n)$ we have
\begin{equation}
\Psi_{\underline{h}} \left(Q,v,k,\gamma \tau , \gamma \underline{z} \right)
=\epsilon (d) (c\tau +d)^{k+r} \exp \left(\frac{\pi icG[\underline{z}] }{c\tau +d}\right) \Psi_{\underline{h}} (Q,v,k,\tau ,\underline{z}) \label{PsiTrans1}
\end{equation}
and
\begin{equation}
 \Psi_{\underline{h}} (Q,v,k,\tau ,\underline{z}+\underline{\lambda} \tau +\underline{\mu})=\exp \left(-\pi i (G[\underline{\lambda}] \tau +2\underline{z}^T G\underline{\lambda})\right) \Psi_{\underline{h}} (Q,v,k,\tau ,\underline{z}). \label{PsiTrans2}
\end{equation}
In particular, $\Psi_{\underline{h}}(Q,v,k,\tau ,\underline{z})$ is a Jacobi form of weight $k+r$, index $G/2$, and character $\epsilon$ on the subgroup $J^n_0 (N)$.
\end{theorem}
 
 To obtain proofs for the previous theorems we consider functions of the form
 \begin{equation}
 \Theta_{\underline{h}} (Q,v,\tau ,\underline{z},X):=\sum_{t\geq 0} \frac{2^{t/2}\theta_{\underline{h}}(Q,v,t,\tau ,\underline{z})}{t!}(2\pi i X)^t , \label{JacobiLikeTheta}
 \end{equation}
 which are similar to Jacobi-like forms. The above theorems follow from manipulation of the following result, which is proved in Section \ref{ProofJacobiLikeTheorem}.

\begin{theorem} \label{JacobiLikeTheorem}
 Suppose $B(v,h_j)=0$ for all $1\leq j\leq n$. Then for any $\gamma =\left( \begin{smallmatrix} a& b\\c &d \end{smallmatrix} \right) \in \Gamma_0 (N)$ we have
 \begin{align}
&\Theta_{\underline{h}} \left( Q,v,\gamma \tau ,\gamma \underline{z},\frac{X}{c\tau +d} \right) \notag \\
 &=\epsilon (d)e^{\frac{\pi icG[\underline{z}]}{c\tau +d}}(c\tau +d)^r \exp \left(\frac{2\pi i[2Q(v)]cX^2}{c\tau +d}\right) 
\Theta_{\underline{h}} \left(Q,v,\tau ,\underline{z},X \right). \label{b1}
 \end{align}
\end{theorem}

 	\indent The author would like to thank Geoffrey Mason and Olav Richter for many valuable discussions. In particular, Richter pointed out that there are other possible approaches to some of the results proved here. The author would also like to thank the referee for their excellent suggestions and remarks.
 

\section{Preliminaries}\label{Preliminaries}

 Let $\text{Hol}_{\mathbb{H}\times \mathbb{C}^n}$ denote the space of holomorphic functions from $\mathbb{H}\times \mathbb{C}^n$ to $\mathbb{C}$, and $F$ be a real symmetric positive definite $n\times n$ matrix. We say a function $\phi \in \text{Hol}_{\mathbb{H}\times \mathbb{C}^n}$ is a \textit{(holomorphic) Jacobi form of weight $k$, index $F$, and character $\chi$ $(\chi \colon \Gamma' \to \mathbb{C}^*)$} on a subgroup $\Gamma'$ of $\Gamma$ if $\phi$ has an expansion of the form
 \begin{equation}
 \phi (\tau ,\underline{z})= \underset{4\ell -F^{-1}[\underline{\nu}] \geq 0}{\sum_{\underline{\nu}\in \mathbb{Z}^n, \ell \in \mathbb{Q},}} c(\ell ,\underline{\nu}) q^\ell \exp \left( 2\pi i (\underline{z}^T \underline{\nu}) \right), \label{JacobiExpansion1}
 \end{equation}
 where $\ell \geq 0$, $c(\ell ,\underline{\nu})$ are scalars, and for all $\gamma =\left( \begin{smallmatrix} a& b\\c &d \end{smallmatrix} \right) \in \Gamma'$ and $(\underline{\lambda},\underline{\mu}) \in \mathbb{Z}^n \times \mathbb{Z}^n$ we have
 \begin{equation}
 \phi (\gamma \tau ,\gamma \underline{z}) = \chi (\gamma) (c\tau +d)^{k} \exp \left( 2\pi i \frac{cF[\underline{z}]}{c\tau +d} \right) \phi (\tau ,\underline{z}), \label{JacobiFormDefn1}
 \end{equation}
 and
 \begin{equation}
  \phi (\tau ,\underline{z}+\underline{\lambda}\tau +\underline{\mu}) = \exp \left(-2\pi i(\tau F[\underline{\lambda}]+2\underline{z}^T F \underline{\lambda}) \right) \phi (\tau ,\underline{z}). \label{JacobiFormDefn2}
 \end{equation}
 
 A function $\phi \in \text{Hol}_{\mathbb{H}\times \mathbb{C}^n}$ is a \textit{quasi-Jacobi form of weight $k$ and index $F$ on $\Gamma'$} if for each $(\tau ,\underline{z}) \in \mathbb{H}\times \mathbb{C}^n$ and all $\gamma =\left( \begin{smallmatrix} a& b\\c &d \end{smallmatrix} \right)\in \Gamma'$ and $(\underline{\lambda},\underline{\mu})\in \mathbb{Z}^n \times \mathbb{Z}^n$, we have 
 \begin{enumerate}
 \item $(c\tau +d)^{-k} \exp \left(-\frac{2\pi icF[\underline{z}]}{c\tau +d} \right) \phi(\gamma \tau ,\gamma \underline{z}) \in \text{Hol}_{\mathbb{H}\times \mathbb{C}^n} \left[\frac{cz_1}{c\tau +d}, \dots ,\frac{cz_n}{c\tau +d}, \frac{c}{c\tau +d}\right]$ with coefficients dependent only on $\phi$, and
 \item $\exp \left( 2\pi i(\tau F[\underline{\lambda}]+2\underline{z}^T F \underline{\lambda}) \right) \phi (\tau ,\underline{z}+\underline{\lambda}\tau +\underline{\mu}) \in \text{Hol}_{\mathbb{H}\times \mathbb{C}^n}[\lambda_1 ,\dots ,\lambda_n]$ with coefficients dependent only on $\phi$.
 \end{enumerate}
  In other words, there are natural numbers $s_1 ,\dots ,s_n ,t$, and holomorphic functions $S_{i_1 ,\dots ,i_n ,j}(\phi)$ and $T_{i_1 ,\dots ,i_n} (\phi)$ on $\mathbb{H} \times \mathbb{C}^n$ determined only by $\phi$ such that
 \begin{align}
 &(c\tau +d)^{-k} \exp \left(-\frac{2\pi icF[\underline{z}]}{c\tau +d} \right) \phi \left( \gamma \tau, \gamma \underline{z} \right)\nonumber \\
 &=\underset{j\leq t}{\sum_{i_1 \leq s_1,\dots ,i_n \leq s_n}} S_{i_1 ,\dots ,i_n,j} (\phi)(\tau ,\underline{z})\left(\frac{cz_1}{c\tau +d}\right)^{i_1}\cdots \left(\frac{cz_n}{c\tau +d}\right)^{i_n} \left( \frac{c}{c\tau +d}\right)^j ,\label{QuasiDefn1}
 \end{align}
 and
 \begin{align}
 & \exp \left( 2\pi i(\tau F[\underline{\lambda}]+2\underline{z}^T F \underline{\lambda}) \right) \phi (\tau ,\underline{z}+\underline{\lambda}\tau +\underline{\mu})\notag \\
 &\hspace{20mm}=\sum_{i_1 \leq s_1 ,\dots ,i_n \leq s_n} T_{i_1 ,\dots ,i_n} (\phi)(\tau ,\underline{z})\lambda_1^{i_1}\cdots \lambda_n^{i_n} . \label{QuasiDefn2}
 \end{align}
 If $\phi \not =0$, we may take $S_{s_1 ,\dots ,s_n ,t}(\phi) \not =0$ and $T_{s_1 ,\dots ,s_n} (\phi) \not =0$ and say that $\phi$ is a quasi-Jacobi form of \textit{depth} $(s_1 ,\dots ,s_n ,t)$.\\ 
 	\indent Direct calculation shows the space of quasi-Jacobi forms is invariant under applications of $\frac{d}{dz_j}$, $\frac{d}{d\tau}$, and $E_2 (\tau)$. In particular, $\frac{d}{dz_j}$ and $E_2 (\tau)$ applied to a quasi-Jacobi form of weight $k$ increases the weight to $k+1$ and $k+2$, respectively.
 
\section{Proofs of Theorems \ref{theorem1} and \ref{theorem2}} \label{ProofTheoremsMajor}

 In this section we assume Theorem \ref{JacobiLikeTheorem} and use this to prove Theorems \ref{theorem1} and \ref{theorem2}. We begin by proving the $\Gamma_0 (N)$ transformations (\ref{ThetaTrans1}) and (\ref{PsiTrans1}). Take $Q(v)\not =0$ and consider the function
 \[
 \widehat{E}_2 \left(Q,v,\tau ,X\right):= \exp \left( 2Q(v)E_2 (\tau) (2\pi i X)^2\right) .
 \]
 Using the transformation $(c\tau +d)^{-2}E_2 (\gamma \tau)=E_2 (\tau) -\frac{c}{2\pi i (c\tau +d)}$, we note
 \begin{equation}
 \widehat{E}_2 \left(Q,v,\gamma \tau ,\frac{X}{c\tau +d}\right)=\widehat{E}_2 \left(Q,v,\tau ,X \right)\exp \left(-\frac{2\pi i c[2Q(v)]X^2}{c\tau +d}\right). \label{49ers1}
 \end{equation}
 Combining (\ref{49ers1}) and Theorem \ref{JacobiLikeTheorem} gives
 \begin{align}
 \widehat{E}_2 & \left(Q,v,\gamma \tau ,\frac{X}{c\tau +d}\right)\Theta_{\underline{h}} \left( Q,v,\gamma \tau ,\gamma \underline{z},\frac{X}{c\tau +d} \right) \notag \\
 &= \widehat{E}_2 \left(Q,v,\tau ,X \right)\epsilon (d)e^{ \frac{\pi icG[\underline{z}]}{c\tau +d}}(c\tau +d)^r \Theta_{\underline{h}} \left(Q,v,\tau ,\underline{z},X \right). \label{bbb2}
 \end{align}
 Expanding the left hand side of (\ref{bbb2}), we find
 \begin{align}
 &\widehat{E}_2 \left(Q,v,\gamma \tau ,\frac{X}{c\tau +d}\right)\Theta_{\underline{h}} \left( Q,v,\gamma \tau ,\gamma \underline{z},\frac{X}{c\tau +d} \right)\notag \\
 &=\sum_{\ell \geq 0} \frac{1}{\ell !} \left(2Q(v)E_2 (\gamma \tau)\right)^\ell \left(\frac{2\pi iX}{c\tau +d}\right)^{2\ell} \sum_{t\geq 0} \frac{2^{t/2}}{t!} \theta_{\underline{h}}(Q,v,t,\gamma \tau ,\gamma \underline{z}) \left(\frac{2\pi iX}{c\tau +d}\right)^t \notag \\
 &=\sum_{\ell \geq 0} \sum_{t\geq 0} \frac{2^{t/2}}{\ell !t!} \left(2Q(v)E_2 (\gamma \tau)\right)^\ell \theta_{\underline{h}}(Q,v,t,\gamma \tau ,\gamma \underline{z}) (c\tau +d)^{-(2\ell+t)} (2\pi iX)^{2\ell +t} \notag \\
 &=\sum_{k\geq 0} \sum_{\ell =0}^{\lfloor k/2 \rfloor} \frac{2^{k/2}}{k!} \delta (\ell ,k) \left(2Q(v)E_2 (\gamma \tau)\right)^\ell \theta_{\underline{h}}(Q,v,k-2\ell ,\gamma \tau ,\gamma \underline{z}) (c\tau +d)^{-k} (2\pi iX)^k, \label{lefthand}
 \end{align}
 where we set $k=2\ell +t$ and use the fact $2^{\frac{1}{2}(k-2\ell)}/ \ell !(k-2\ell)!=2^{k/2}\delta (\ell ,k)/k!$. Expanding the right hand side of (\ref{bbb2}) shows
 \begin{align}
  &\widehat{E}_2 \left(Q,v,\tau ,X \right)\epsilon (d)e^{ \frac{\pi icG[\underline{z}]}{c\tau +d}}(c\tau +d)^r \Theta_{\underline{h}} \left(Q,v,\tau ,\underline{z},X \right)\notag  \\
  &=\epsilon (d)(c\tau +d)^r e^{ \frac{\pi icG[\underline{z}]}{c\tau +d}} \sum_{\ell \geq 0} \sum_{t\geq 0} \frac{2^{t/2}}{\ell !t!} \left(2Q(v)E_2(\tau)\right)^\ell  \theta_{\underline{h}}(Q,v,t,\tau ,\underline{z}) (2\pi i X)^{2\ell +t}\notag  \\
  &=\epsilon (d)(c\tau +d)^r e^{ \frac{\pi icG[\underline{z}]}{c\tau +d}} \sum_{k\geq 0} \sum_{\ell =0}^{\lfloor k/2 \rfloor} \frac{2^{k/2}}{k!} \delta (\ell ,k) \left(2Q(v)E_2 (\tau)\right)^\ell \notag \\
  &\hspace{20mm}\cdot \theta_{\underline{h}}(Q,v,k-2\ell ,\tau ,\underline{z}) (2\pi iX)^k , \label{righthand}
 \end{align}
 where we again set $k=2\ell +t$. Using (\ref{bbb2}) to combine (\ref{lefthand}) and (\ref{righthand}) and then comparing the coefficients of $X^k$, we obtain
\begin{align*}
  &\frac{2^{k/2}}{k!} \sum_{\ell =0}^{\lfloor k/2 \rfloor} \delta (\ell ,k) \left(2Q(v)E_2 (\gamma \tau)\right)^\ell \theta_{\underline{h}}(Q,v,k-2\ell ,\gamma \tau ,\gamma \underline{z}) \\
  &\hspace{5mm}=\epsilon (d)(c\tau +d)^{r+k} e^{ \frac{\pi icG[\underline{z}]}{c\tau +d}} \frac{2^{k/2}}{k!} \sum_{\ell =0}^{\lfloor k/2 \rfloor}  \delta (\ell ,k) \left(2Q(v)E_2 (\tau)\right)^\ell \theta_{\underline{h}}(Q,v,k-2\ell ,\tau ,\underline{z}).
\end{align*}
That is,
\[
\Psi_{\underline{h}}(Q,v,k,\gamma \tau ,\gamma \underline{z})=\epsilon (d)(c\tau +d)^{r+k} e^{ \frac{\pi icG[\underline{z}]}{c\tau +d}} \Psi_{\underline{h}}(Q,v,k,\tau ,\underline{z}),
\]
as desired. This establishes (\ref{PsiTrans1}). Taking $Q(v)=0$ gives (\ref{ThetaTrans1}). \\

 We will now prove the $\mathbb{Z}^n \times \mathbb{Z}^n$ transformations (\ref{ThetaTrans2}) and (\ref{PsiTrans2}). Since $\mu_j B(m,h_j) \in \mathbb{Z}$ for each $j$, we have $\theta_{\underline{h}}(Q,v,k,\tau ,\underline{z}+\underline{\lambda}\tau +\underline{\mu})=\theta_{\underline{h}}(Q,v,k,\tau ,\underline{z}+\underline{\lambda}\tau)$. Using that
 \[
 Q(m+\lambda_1 h_1 +\cdots +\lambda_n h_n) =Q(m) + \sum_{j=1}^n \lambda_j B(m,h_j ) +\frac{1}{2}G[\underline{\lambda}],
 \]
 we find
 \begin{align}
 &\theta_{\underline{h}}(Q,v,k,\tau ,\underline{z}+\underline{\lambda}\tau +\underline{\mu})  \notag \\
 &=e^{-\pi i \tau G[\underline{\lambda}]} \sum_{m\in \mathbb{Z}^f} B(v,m)^k e^{2\pi i \tau Q\left( m+\sum_{j=1}^n \lambda_j h_j\right)} e^{2\pi i \sum_{j=1}^n z_j B(m,h_j)} \notag \\
 &=e^{-\pi i (G[\underline{\lambda}]\tau -2\underline{z}^T G\underline{h})} \sum_{m\in \mathbb{Z}^f} B\left( v,m-\sum_{j=1}^n \lambda_j h_j\right)^k e^{2\pi i Q(m)}e^{2\pi i \sum_{j=1}^n z_j B(m,h_j)} , \notag
 \end{align}
 where we replaced $m$ with $m-\sum_{j=1}^n \lambda_j h_j$ in the last equality. Using the assumption $B(v,h_j)=0$ for $1\leq j\leq n$ establishes (\ref{ThetaTrans2}). Equation (\ref{PsiTrans2}) follows immediately from (\ref{ThetaTrans2}). The proof of Theorem \ref{theorem2} is complete.\\
 
 We now turn our attention to proving (\ref{thm1a}) and (\ref{Expansion}), and therefore no longer assume $B(v,h_j)=0$ for all $j$. First we establish the following lemma.
 
\begin{lemma} \label{LemmaQuasi}
 For any $v\in \mathbb{C}^f$ there exists an element $u \in \mathbb{C}^f$ satisfying $B(u,h_j)=0$ for $1\leq j\leq n$ such that 
 \[
 \theta_{\underline{h}}(Q,v,k,\tau ,\underline{z})
 =\left( \sum_{p_1+\cdots +p_n +k_1  =k} \binom{k}{p_1,\dots ,p_n,k_1} \prod_{i=1}^n \alpha_i^{p_i} \partial_{z_i}^{p_i}\right) \theta_{\underline{h}}\left(Q, u ,k_1 ,\tau ,\underline{z}\right).
 \]
 Here $\underset{p_1+\cdots +p_n +k_1  =k}{\sum}$ denotes summing over the positive integers $p_1, \dots ,p_n$, and $k_1$ which sum to $k$, $\alpha_i$ are scalars, $\partial_{z_i} =\frac{1}{2\pi i}\frac{d}{d z_i}$, and $\binom{k}{p_1,\dots ,p_n,k_1}$ are are the multinomial coefficients $\binom{k}{p_1,\dots ,p_n,k_1} =\frac{k!}{p_1! \cdots p_n! k_1!}$.
\end{lemma}

\paragraph{Proof}
 Extend the set $\{h_1 ,\dots ,h_n \}$ into a basis $\{h_1 ,\dots ,h_n ,u_{n+1} ,\dots ,u_f \}$ of $\mathbb{C}^f$ such that $B(h_i ,u_j)=0$ for all $i,j$. Then there are scalars $\alpha_i ,\beta_j \in \mathbb{C}$ such that $v=\sum_{i} \alpha_i h_i + \sum_j \beta_j u_j$. Set $u=\sum_j \beta_j u_j$. We have
 \[
 \theta_{\underline{h}}(Q,v,k,\tau ,\underline{z})
  =\sum_{m\in \mathbb{Z}^f} B\left(\sum_{i} \alpha_i h_i +u ,m \right)^k q^{Q(m)} \zeta_1^{B(m,h_1)} \cdots \zeta_n^{B(m,h_n)}.
 \]
 \pagebreak
 The lemma follows by use of the multinomial theorem and replacing each $B(h_i ,m)$ with $\partial_{z_i}$.
 \vspace{-7mm} \begin{flushright}$\Box$ \end{flushright}
 
 We now establish the convergence of the functions $\theta_{\underline{h}}(Q,v,k,\tau ,\underline{z})$. If we again extend the set $\{h_1, \dots ,h_n \}$ into a basis $\{h_1 ,\dots ,h_n ,u_{n+1}, \dots ,u_f \}$ for $\mathbb{C}^f$ and write $v=\sum_{i} \alpha_i h_i + \sum_j \alpha_j u_j$, we may consider functions of the form
 \[
 \theta_{\underline{h},\underline{u}}(Q,v,k,\tau ,\underline{z})=\sum_{m\in \mathbb{Z}^f} B(v,m)^k q^{Q(m)} \zeta_1^{B(m,h_1)} \cdots \zeta_n^{B(m,h_n)} \zeta_{n+1}^{B(m,u_{n+1})} \cdots \zeta_f^{B(m,u_f)}.
 \]
  By the previous lemma we have
  \[
 \theta_{\underline{h},\underline{u}}(Q,v,k,\tau ,\underline{z})
 =\left( \sum_{p_1+\cdots +p_f=k} \binom{k}{p_1,\dots ,p_f} \prod_{i=1}^f \alpha_i^{p_i} \partial_{z_i}^{p_i}\right) \theta_{\underline{h},\underline{u}}\left(Q ,0 ,\tau ,\underline{z}\right),
 \]
 where $\theta_{\underline{h},\underline{u}}\left(Q ,0,\tau ,\underline{z}\right)=\sum_{m\in \mathbb{Z}^f} q^{Q(m)} \zeta_1^{B(m,h_1)} \cdots \zeta_n^{B(m,h_n)}\zeta_{n+1}^{B(m,u_{n+1})} \cdots \zeta_f^{B(m,u_f)}$ converges for $\tau \in \mathbb{H}$ and each $z_i \in \mathbb{C}$ when we fix the remaining $z_j$, $j\not =i$. This last statement is dependent on the fact $Q$ is positive definite, so that $\vert q^{Q(m)}\vert <1$ for all nonzero $m\in \mathbb{Z}^f$. Therefore, by Hartogs' Theorem the function converges for all $(\tau ,\underline{z})$ on $\mathbb{H}\times \mathbb{C}^f$, and it follows that $\theta_{\underline{h},\underline{u}}(Q,v,k,\tau ,\underline{z})$ converges on $\mathbb{H} \times \mathbb{C}^f$. Setting $z_{n+1}=\cdots =z_f =0$ shows that $\theta_{\underline{h}}(Q,v,k,\tau ,\underline{z})$ is convergent on $\mathbb{H}\times \mathbb{C}^n$.\\
 
 We turn to proving the inequality
 \[
 4Q(m) -(B(m,h_1), \dots ,B(m,h_n))(G/2)^{-1} (B(m,h_1), \dots ,B(m,h_n))^T \geq 0.
 \]
 Setting $\beta =(B(m,h_1), \dots ,B(m,h_n))$, we rewrite this expression as
 \begin{equation}
 \beta G^{-1} \beta^T \leq B(m,m). \label{BilinearFormEquality}
 \end{equation}
 Since $G$ is a real symmetric matrix it has the decomposition $G=QDQ^{-1} =QDQ^T$, where $D$ is a diagonal matrix and $Q=(q_{ij})$ is orthogonal. Therefore, $G^{-1} =QD^{-1} Q^T$ and we rewrite (\ref{BilinearFormEquality}) as $(\beta Q)D^{-1} (\beta Q)^T \leq B(m,m)$. Setting $v_i =\sum_j q_{ji}h_j$, this becomes
 \begin{equation}
  (B(m,v_1) ,\dots ,B(m,v_n))D^{-1} (B(m,v_1) ,\dots ,B(m,v_n))^T \leq B(m,m). \label{BilinearFormEquality3}
 \end{equation}
 So long as $D= (B(v_i ,v_j))$, establishing (\ref{BilinearFormEquality}) is equivalent to proving (\ref{BilinearFormEquality3}). However, this is the case since 
 \[
 B(v_i,v_j)=B(v_j,v_i) = \left( \sum_r q_{rj}h_r ,\sum_s q_{si}h_s \right)=\sum_{r,s} q_{rj}B(h_r,h_s)q_{si}
 \]
  is the $(i,j)$-th component of the matrix $D=Q^T GQ$. It is therefore sufficient to verify (\ref{BilinearFormEquality3}).\\
  \indent First we assume that $h_1 ,\dots ,h_n$ span $\mathbb{Z}^f$. Then $m=\sum_j \lambda_j v_j$ for some $\lambda_j$. In this case we have
  \begin{align*}
  &(B(m,v_1) ,\dots ,B(m,v_n))D^{-1} (B(m,v_1) ,\dots ,B(m,v_n))^T \\
  &\hspace{10mm}=\sum_{i,j} \lambda_j^2 B(v_j ,v_i)^2 B(v_i ,v_i)^{-1} =\sum_{i} \lambda_i^2 B(v_i ,v_i) =B(m,m)
  \end{align*}
  since $B(v_j ,v_i)=0$ for $i\not =j$, being the off-diagonal components of the diagonal matrix $D$. This establishes (\ref{BilinearFormEquality3}) when $h_1 ,\dots ,h_n$ span $\mathbb{Z}^f$.\\
  \indent Suppose next that $h_1 ,\dots ,h_n$ do not span $\mathbb{Z}^f$. Since they are linear independent they form an $\mathbb{R}$-basis for $\mathbb{R}^n$. We consider the orthogonal semi-direct product $\mathbb{R}^f =\mathbb{R}^n \perp \mathbb{R}^{f-n}$. Let $h_{n+1} ,\dots ,h_f$ be a basis for $\mathbb{R}^{f-n}$ and write $m=m' +m''$, where $m' \in \mathbb{R}^n$ and $m'' \in \mathbb{R}^{f-n}$. Note that $B(m,v_i)=B(m' +m'',v_i)=B(m' ,v_i)$ and $B(m' ,m'')=0$, so that establishing (\ref{BilinearFormEquality3}) is equivalent to proving
  \begin{align*}
  (B(m',v_1) ,\dots ,B(m',v_n))&D^{-1} (B(m',v_1) ,\dots ,B(m',v_n))^T \\
  &\leq B(m',m')+B(m'' ,m'').
  \end{align*}
  However, we find
  \[
  (B(m',v_1) ,\dots ,B(m',v_n))D^{-1} (B(m',v_1) ,\dots ,B(m',v_n))^T = B(m',m')
  \]
  by our arguments above and the fact $B(m'',m'')\geq 0$. This proves (\ref{BilinearFormEquality}) for all linearly independent $h_1 ,\dots ,h_n$.\\
  \indent Note that it is the expansion (\ref{Expansion}) which allows us to claim the functions $\Psi_{\underline{h}}(Q,v,k,\tau ,\underline{z})$ in Theorem \ref{theorem2} are Jacobi forms. In particular, since $E_2 (\tau)$ has a Fourier expansion with positive powers of $q$ and the Fourier expansion of $\theta_{\underline{h}}(Q,v,k-2t,\tau ,\underline{z})$ satisfies (\ref{BilinearFormEquality}), it must be that $\Psi_{\underline{h}}(Q,v,k,\tau ,\underline{z})$ satisfies (\ref{JacobiExpansion1}).\\
  
 \indent It remains to prove (\ref{thm1a}). As mentioned before, $\partial_{z_i}$ maps quasi-Jacobi forms of weight $k$ to weight $k+1$. In the case $Q(v)=0$, (\ref{Expansion}), (\ref{ThetaTrans1}), and (\ref{ThetaTrans2}), along with Lemma \ref{LemmaQuasi}, establish that $\theta_{\underline{h}}(Q,v,k,\tau ,\underline{z})$ is a quasi-Jacobi form of weight $k+r$ and index $G/2$.\\
 \indent Assume now that $Q(v)\not =0$. Recalling Lemma \ref{LemmaQuasi} (and its notation) we have
 \begin{align}
 &\theta_{\underline{h}}(Q,v,k,\tau ,\underline{z}) \label{vtou} \\
 & =\left( \sum_{p_1+\cdots +p_n +k_1  =k} \binom{k}{p_1,\dots ,p_n,k_1} \prod_{i=1}^f \alpha_i^{p_i} \partial_{z_i}^{p_i}\right) \theta_{\underline{h}}(Q,u,k_1 ,\tau ,\underline{z}). \notag
 \end{align}
 \indent We will prove by induction on $k_1$ that $\theta_{\underline{h}}(Q,u,k_1,\tau ,\underline{z})$ is a quasi-Jacobi form of weight $k_1 +r$ and index $G/2$. Clearly $\theta_{\underline{h}}(Q,u,0,\tau ,\underline{z})$ is a Jacobi form of weight $r$ and index $G/2$. Suppose that for each $\ell <k_1$, $\theta_{\underline{h}}(Q,u,\ell ,\tau ,\underline{z})$ is a quasi-Jacobi form of weight $\ell +r$ and index $G/2$. Since $\delta (0,k_1)=1$ and $B(u,h_j)=0$ for each $j$, Theorem \ref{theorem2} gives
 \begin{align*}
\Psi_{\underline{h}}&(Q,u,k_1 ,\tau ,\underline{z}) \\
 &=\theta_{\underline{h}}(Q,u,k_1,\tau ,\underline{z}) +\sum_{t=1}^{\lfloor k_1 /2 \rfloor} \delta (t,k_1) (2Q(u)E_2 (\tau))^{t} \theta_{\underline{h}}(Q,u,k_1 -2t,\tau ,\underline{z})
 \end{align*}
 is a Jacobi form of weight $k_1 +r$. Using our induction hypothesis we find $\theta_{\underline{h}}(Q,u,k_1 ,\tau ,\underline{z})$ is a quasi-Jacobi form of weight $k_1 +r$ and index $G/2$. Combining this with (\ref{vtou}) completes the proof of (\ref{thm1a}).
\section{Proof of Theorem \ref{JacobiLikeTheorem}}\label{ProofJacobiLikeTheorem}
 
 Let $A$ be a matrix of level $N$ such that $2Q(x)=x^T Ax$ and $B(x,y)=x^T Ay$. Take $p$ so that $Ap\equiv 0 \mod N$. For $\ell^T =(\ell_1 ,\dots ,\ell_f) \in \mathbb{C}^f$ we set
\begin{align}
 &\theta (A,p,\ell ,k,\tau ,\underline{z})\label{DM-3.9} \\
 &=\frac{1}{N^k}\underset{m\equiv p(N)}{\sum_{m\in \mathbb{Z}^f}} (\ell^T Am)^{k} e^{2\pi i \left(\tau Q(m)/N^2 + z_1 B(m,h_1)/N+\cdots + z_n B(m,h_n )/N \right)} . \notag
\end{align}

 The following theorem is analogous to Theorem $3.4$ in \cite{DM-Theta} but includes the complex variables $\underline{z}$\ and vectors $h_1 ,\dots ,h_n$. It is important for the proof of Theorem \ref{JacobiLikeTheorem}.

\begin{theorem} \label{DM-Theorem-3.4}
 Let the notation be as before and assume $\ell^T A h_j =0$ for all $1\leq j \leq n$ (that is, $B(\ell ,h_j)=0$). If $\gamma =\left( \begin{smallmatrix} a& b\\c &d \end{smallmatrix} \right) \in \Gamma_0 (N)$ and $d>0$, then
 \begin{align}
 &e^{ \frac{-\pi ic G[\underline{z}]}{c\tau +d}}(c\tau +d)^{-(r+k)} \theta \left( A,p,\ell ,k, \gamma \tau ,\gamma \underline{z} \right) \notag \\
 &\hspace{1mm}=\epsilon (d) \exp (2\pi i Q(p)ab/N^2) \sum_{j=0}^{\lfloor k/2 \rfloor} \left( \frac{Q(\ell)c}{\pi i (c\tau +d)}\right)^j \delta (j,k) \theta (A, bp,\ell ,k-2j,\tau ,\underline{z}). \label{DM-3.13}
 \end{align}
 In particular, if we take $p=0$, then
 \begin{align}
 e^{\frac{-\pi ic G[\underline{z}]}{c\tau +d}}& (c\tau +d)^{-(r+k)} \theta \left( A,\ell ,k, \gamma \tau ,\gamma \underline{z} \right)\notag \\
 & =  \epsilon (d) \sum_{j=0}^{\lfloor k/2 \rfloor} \left( \frac{Q(\ell)c}{\pi i (c\tau +d)}\right)^j \delta (j,k) \theta (A,\ell ,k-2j,\tau ,\underline{z}). \label{DM-3.14}
 \end{align}
\end{theorem}

\paragraph{Proof}
 The ideas and techniques parallel those in \cite{DM-Theta}, which in turn are based on \cite{Schoeneberg-Elliptic}. We begin with a function of the form $\Theta_{\underline{h}} (Q,x):=\sum_{m\in \mathbb{Z}^f} e^{-2 Q(m+x)}e^{-2 B(m+x,h)}$ and consider its Fourier coefficients. The analysis is similar to that of \cite{DM-Theta} and \cite{Schoeneberg-Elliptic} with occasional changes and we omit the details. One key difference, however, is that we must replace $h$ with $z_1h_1 +\cdots +z_n h_n$ during the proof. The reader may also consult \cite{KrauelThesis} for a detailed proof of this theorem. \vspace{-7mm} \begin{flushright}$\Box$ \end{flushright}

 \indent Define the functions
 \[
 \Theta_{\underline{h}}^{\text{even}}(Q,v,\tau ,\underline{z},X):=\sum_{t\geq 0} \frac{2^t \theta_{\underline{h}}(Q,v,2t,\tau ,\underline{z})}{(2t)!}(2\pi i X)^{2t} ,
 \]
 and
 \[
  \Theta_{\underline{h}}^{\text{odd}}(Q,v,\tau ,\underline{z},X):=\sum_{t\geq 0} \frac{2^{(t+1/2)} \theta_{\underline{h}}(Q,v,2t+1,\tau ,\underline{z})}{(2t+1)!}(2\pi i X)^{2t+1} .
 \]
 Note that $\Theta_{\underline{h}}(Q,v,\tau ,\underline{z},X)=\Theta_{\underline{h}}^{\text{even}}(Q,v,\tau ,\underline{z},X) + \Theta_{\underline{h}}^{\text{odd}}(Q,v,\tau ,\underline{z},X)$ and $\theta_{\underline{h}}(Q,v,k,\tau ,\underline{z})=\theta (Q,0,v,k,\tau,\underline{z})$. Using Theorem \ref{DM-Theorem-3.4}, we find for functions of the form (\ref{JacobiLikeTheta}) that
\begin{align*}
 \Theta_{\underline{h}} &\left( Q,v,\gamma \tau ,\gamma \underline{z},\frac{X}{c\tau +d} \right)
 =\sum_{t\geq 0} \frac{2^{t/2} \theta_{\underline{h}} (Q,v,t,\gamma \tau ,\gamma \underline{z})}{t!} \left( \frac{2\pi iX}{c\tau +d} \right)^{t} \\
 &=\epsilon (d) e^{\frac{\pi icG[\underline{z}]}{c\tau +d}}(c\tau +d)^r \sum_{t\geq 0} \sum_{j=0}^{\lfloor t/2 \rfloor} \frac{2^{t/2}}{t!} \left( \frac{2Q(v)c}{2\pi i(c\tau +d)}\right)^j \\
 &\hspace{15mm} \cdot \delta (j,t)(c\tau +d)^{t} \theta_{\underline{h}} (Q,v,t-2j,\tau ,\underline{z}) \left( \frac{2\pi iX}{c\tau +d}\right)^{t} \\
 &= \epsilon (d)e^{\frac{\pi icG[\underline{z}]}{c\tau +d}}(c\tau +d)^r \sum_{t\geq 0} \sum_{j=0}^{\lfloor t/2 \rfloor} \frac{2^{\frac{1}{2}(t-2j)}}{(t-2j)!}\\
 &\hspace{15mm}\cdot \theta_{\underline{h}} (Q,v,t-2j,\tau ,\underline{z})\left(2\pi i X \right)^{t-2j} \left( \frac{2Q(v)c}{2\pi i(c\tau +d)}\right)^j \frac{1}{j!}\left(2\pi i X \right)^{2j} \\ 
 &=\epsilon (d) e^{\frac{\pi icG[\underline{z}]}{c\tau +d}}(c\tau +d)^r \sum_{t\geq 0} \sum_{j=0}^{\lfloor t/2 \rfloor} \frac{2^{\frac{1}{2}(t-2j)}}{(t-2j)!} \theta_{\underline{h}} (Q,v,t-2j,\tau ,\underline{z})  \\
 &\hspace{15mm} \cdot \left(2\pi i X \right)^{t-2j} \frac{1}{j!} \left( \frac{2\pi i[2Q(v)]cX^2}{c\tau +d}\right)^j , 
  \end{align*}
 where we used $\frac{2^{t/2}}{t!} \delta (j,t)=\frac{2^{(t-2j)/2}}{(t-2j)!}\frac{1}{j!}$. Continuing the calculation we find
 \begin{align*}
 \Theta_{\underline{h}}&\left( Q,v,\gamma \tau ,\gamma \underline{z},\frac{X}{c\tau +d} \right)\\
 &= \epsilon (d)e^{\frac{\pi icG[\underline{z}]}{c\tau +d}}(c\tau +d)^r \Biggl [
 \sum_{t\geq 0} \sum_{j=0}^{\lfloor t/2 \rfloor} \frac{2^{\frac{1}{2}(2t-2j)}}{(2t-2j)!} \\
 &\hspace{15mm} \cdot \theta_{\underline{h}} (Q,v,2t-2j,\tau ,\underline{z})(2\pi i X)^{2t-2j} \frac{1}{j!} \left( \frac{2\pi i[2Q(v)]cX^2}{c\tau +d}\right)^j \\
 &\hspace{6mm} +\sum_{t\geq 0} \sum_{j=0}^{\lfloor t/2 \rfloor} \frac{2^{\frac{1}{2}(2t-2j+1)}}{(2t-2j+1)!} \theta_{\underline{h}} (Q,v,2t-2j+1,\tau ,\underline{z})\\
 &\hspace{15mm} \cdot (2\pi i X)^{2t-2j+1} \frac{1}{j!} \left( \frac{2\pi i[2Q(v)]cX^2}{c\tau +d}\right)^j \Biggl ]\\
 &=\epsilon (d)e^{\frac{\pi icG[\underline{z}]}{c\tau +d}}(c\tau +d)^r \exp \left(\frac{2\pi i[2Q(v)]cX^2}{c\tau +d} \right) \\
 &\hspace{15mm}\cdot \biggl [ \Theta_{\underline{h}}^{\text{even}}(Q,v,\tau ,\underline{z},X) +\Theta_{\underline{h}}^{\text{odd}}(Q,v,\tau ,\underline{z},X)\biggl ] \\
 &=\epsilon (d)e^{\frac{\pi icG[\underline{z}]}{c\tau +d}}(c\tau +d)^r \exp \left(\frac{2\pi i[2Q(v)]cX^2}{c\tau +d}\right) \Theta_{\underline{h}} \left(Q,v,\tau ,\underline{z},X \right).
\end{align*}
 We have now proved Theorem \ref{JacobiLikeTheorem} for matrices $\gamma_1 =\left( \begin{smallmatrix} a& b\\c &d \end{smallmatrix} \right) \in \Gamma_0 (N)$ with $d>0$. Let $\gamma_2 =-\gamma_1$. Note $\gamma_2 \tau =\gamma_1 \tau$ and $\gamma_2 \underline{z} =-\gamma_1 \underline{z}$, so that we have
\[
\theta_{\underline{h}}(Q,v,k,\gamma_2 \tau ,\gamma_2 \underline{z})=\theta_{\underline{h}}(Q,v,k, \gamma_1 \tau ,-\gamma_1 \underline{z})=(-1)^k \theta_{\underline{h}}(Q,v,k, \gamma_1 \tau ,\gamma_1 \underline{z}).
\]
It follows that
\[
\Theta_{\underline{h}} \left( Q,v,\gamma_2 \tau ,\gamma_2 \underline{z},\frac{X}{-c\tau -d} \right)
=\Theta_{\underline{h}} \left( Q,v,\gamma_1 \tau ,\gamma_1 \underline{z},\frac{X}{c\tau +d} \right).
\]
Using (\ref{b1}) which we have already established for $d>0$, we find
\begin{align}
&\Theta_{\underline{h}} \left( Q,v,\gamma_2 \tau ,\gamma_2 \underline{z},\frac{X}{-c\tau -d} \right) \notag \\
&\hspace{1mm}=\epsilon (-d)e^{ \frac{\pi icG[\underline{z}]}{c\tau +d}}(-c\tau -d)^r\exp \left(\frac{2\pi i[2Q(v)]cX^2}{c\tau +d}\right) \Theta_{\underline{h}} \left(Q,v,\tau ,\underline{z},X \right). \label{finaleq}
\end{align}
However, (\ref{finaleq}) is also obtained if we replace $c$ and $d$ on the right hand side of (\ref{b1}) with $-c$ and $-d$, respectively. It follows that (\ref{b1}) also holds for $\gamma_2$, and so for all $\gamma \in \Gamma_0 (N)$. The proof of Theorem \ref{JacobiLikeTheorem} is now complete. \vspace{-7mm} \begin{flushright}$\Box$ \end{flushright}



\begin{thebibliography}{0}

\bibitem{AM-BehaviorI}
A.~Andriaov and G.~Malolekin,
\newblock Behavior of theta series of degree $N$ under modular substitutions,
\newblock in {\em Math. USSR-Izv}, \textbf{9}(2) (1975) 227--241.

\bibitem{DM-Theta}
C.~Dong and G.~Mason,
\newblock Transformation laws for theta functions,
\newblock in {\em Proceedings on {M}oonshine and related topics ({M}ontr\'eal,
  {QC}, 1999)}, ed.~J. McKay and A. Sebbar, CRM Proc. Lecture Notes, Vol. 30, (Amer. Math. Soc., 2001) 15--26. 

\bibitem{DMN-quasi}
C.~Dong, G.~Mason, and K.~Nagatomo,
\newblock Quasi-modular forms and trace functions associated to free boson and
  lattice vertex operator algebras,
\newblock in {\em Internat. Math. Res. Notices}, (8) (2001) 409--427.

\bibitem{EZ}
M.~Eichler and D.~Zagier,
\newblock {\em The theory of {J}acobi forms}, Vol. 55 of {\em Progress in
  Mathematics},
\newblock Birkh\"auser Boston Inc., Boston, MA, (1985).

\bibitem{Hecke-Positiven}
E.~Hecke,
\newblock Analytische {A}rithmetik der positiven quadratischen {F}ormen,
\newblock in {\em Danske Vid. Selsk. Math.-Fys. Medd.}, \textbf{17}(12) (1940) 134.

\bibitem{KrauelThesis}
M.~Krauel,
\newblock Vertex operator algebras and {J}acobi forms, ({P}h{D}
  {D}issertation),
\newblock in {\em ProQuest Dissertations and Theses (1039264015). ISBN:
  9781267533739}, (2012).

\bibitem{KrauelMasonI}
M.~Krauel and G.~Mason,
\newblock Vertex operator algebras and weak {J}acobi forms,
\newblock in {\em Int. Journ. Math.}, \textbf{23}(6) (2012).
%

\bibitem{Ogg-Modular}
A.~Ogg,
\newblock {\em Modular forms and {D}irichlet series},
\newblock W. A. Benjamin, Inc., New York-Amsterdam, (1969).

\bibitem{Richter-OnTransformations}
O.~Richter,
\newblock On transformation laws for theta functions,
\newblock in {\em Rocky Mountain J. Math.}, \textbf{34}(4) (2004) 1473--1481.

\bibitem{Schoeneberg-Das}
B.~Schoeneberg,
\newblock Das {V}erhalten von mehrfachen {T}hetareihen bei
  {M}odulsubstitutionen,
\newblock in {\em Math. Ann.}, \textbf{116}(1) (1936) 511--523.

\bibitem{Schoeneberg-Elliptic}
B.~Schoeneberg,
\newblock {\em Elliptic modular functions: an introduction},
\newblock Springer-Verlag, New York, (1974).
\newblock Translated from the German by J. R. Smart and E. A. Schwandt, Die
  Grundlehren der mathematischen Wissenschaften, Band 203.
  
\bibitem{Ziegler-JacobiForms}
C.~Ziegler,
\newblock Jacobi forms of higher degree,
\newblock in {\em Abh. Math. Sem. Univ. Hamburg}, (59) (1989) 191--224.

\end{thebibliography}
\end{document}